\documentclass[a4paper,11pt]{article}
  \usepackage[dvips]{graphicx}
\usepackage{amsmath}
\usepackage{url,mstorti}
\usepackage{hyperref,nameref}
\usepackage{natbib}

\newcommand{\MMw}{{\mathrm{MM}_\mathrm{w}}}
\newcommand{\MMa}{{\mathrm{MM}_\mathrm{a}}}
\newcommand{\Cp}{{C_p}}

\newcommand{\Cpv}{{C_{pv}}}
\newcommand{\vg}{\mathbf{v}\G}
\newcommand{\Vg}{\mathbf{V}\G}

\newcommand{\Ql}{Q_l}
\newcommand{\Dv}{D_v}
\newcommand{\Da}{D_a}
\newcommand{\Pv}{P_v}
\newcommand{\Psat}{P_\mathrm{sat}}

\newcommand{\Kg}{K_g}
\newcommand{\KKg}{{!K_g}}
\newcommand{\HH}{\mathrm{H}}
\newcommand{\HR}{\mathrm{H\!R}}
\newcommand{\Dt}{\mathrm{\Delta t}}
\newcommand{\ave}[1]{\left<#1\right>}

\newcommand{\atm}{{\mathrm{atm}}}

\newcommand{\keff}{k_{\mathrm{eff}}}
\newcommand{\Tref}{T_{\mathrm{ref}}}
\newcommand{\effe}{{\mathrm{eff}}}
\newcommand{\Tplaten}{T_{\mathrm{platen}}}
\newcommand{\be}{\begin{equation}}
\newcommand{\ee}{\end{equation}}
\newcommand{\bea}{\begin{eqnarray}}
\newcommand{\eea}{\end{eqnarray}}

\newcommand{\ben}{\begin{enumerate}}
\newcommand{\een}{\end{enumerate}}
\newcommand{\bit}{\begin{itemize}}
\newcommand{\eit}{\end{itemize}}
\newcommand{\ba}{\begin{array}}
\newcommand{\ea}{\end{array}}

\newcommand{\bc}{\begin{center}}
\newcommand{\ec}{\end{center}}

\newcommand{\eg}{{\epsilon\G}}
\newcommand{\rhog}{{\rho\G}}
\newcommand{\G}{{_g}}
\newcommand{\So}{{_s}}
\newcommand{\es}{{\epsilon\So}}
\newcommand{\rhos}{{\rho\So}}
\newcommand{\Li}{{_l}}
\newcommand{\el}{{\epsilon\Li}}
\newcommand{\rhol}{{\rho\Li}}
\newcommand{\hl}{{h\Li}}

\newcommand{\hg}{{h\G}}
\begin{document}
\sloppy

\title{Hot-pressing process modeling for medium density \\
        fiberboard (MDF)}
\author{
Norberto Nigro and Mario Storti\\ 
Centro Internacional de M\'etodos Computacionales en Ingenier\'\i{}a\\
\url{http://www.cimec.com.ar}, \url{<mstorti@intec.unl.edu.ar>}}
\date{\today}
\maketitle

\abstract{ 
In this paper we present a numerical solution for the
mathematical modeling of the hot-pressing process applied to medium
density fiberboard. The model is based in the work of
\citet{humphrey82,humphrey89} and \citet{carvalho98}, with some modifications and
extensions in order to take into account mainly the convective effects
on the phase change term and also a conservative numerical treatment
of the resulting system of partial differential equations. 
}

\tableofcontents

\section*{Classification}

\begin{itemize}
\item[\tt 65M60] \small Finite elements, Rayleigh-Ritz and Galerkin methods, finite methods 
\item[\tt 82C70] Transport processes
\item[\tt 76S05] Flows in porous media; filtration; seepage 
\item[\tt 76T30] Three or more component flows
\end{itemize}

\section{Hot-pressing mathematical model}

\begin{figure*}[htb]
\centerline{\includegraphics{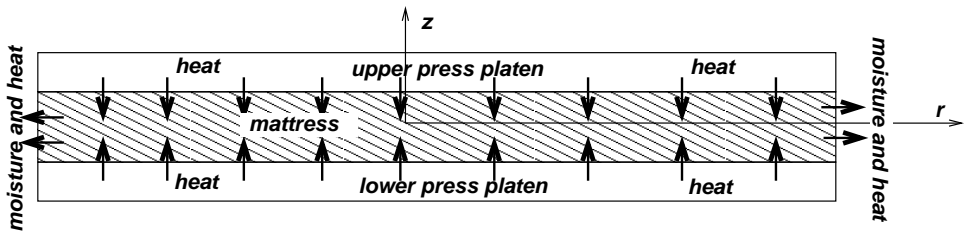}}
\caption{Hot pressing process description}
\label{fg:hotpress}
\end{figure*}

Hot pressing is the process in which a mattress composed of wood
fibers and resin is cured by applying heat and pressure in a press
(see figure~\ref{fg:hotpress}). Continuum and batch presses do exist, and
one of the main issues in reducing the cost of the final product is
to reduce the press cycle time. In order to improve the heat transfer
between the press platens and the inner layers, some amount of water
is added to the mat. Another issue is to adjust the parameters and the
temperature history of the cycle in order to obtain a given density
profile in the board. Normally, it is desirable to have lower
densities at the center of the board in order to increase the
mechanical rigidity for a given total mass per unit area. Predicting
the influence of these parameters, namely water content, press cycle
duration and history (pressure and temperature) is one of the main
concerns of numerical models. 

Many numerical models have been reported to help in predicting the
influence of the process parameters in the final product. Among the
most complete, we can find the finite difference 2D (axisymmetrical)
model presented by \citet{humphrey82} and, more recently the 3D model
of \citet{carvalho98}. Both  consider conduction, phase
change of water from the adsorbed to the vapor state and
convection. The stress development and the determination of
the density profile are not included in these models. 

In this paper we present a numerical model which includes all these
features and makes some correction to the energy balance equation, as
presented in \citet{carvalho98}. The model is based on the finite
element model, so that it allows for a more versatile definition  of 
geometry, dimensionality and,  eventually, coupling with other
packages. It also will allow the use of adaptive refinement, which may
be an important issue at the lateral borders, where the hot steam
flows from the board to the ambient. However, this issue is not
considered in this paper. 

\subsection{Multiphase model}

In order to avoid modeling the material down to the scale of the
microstructure (the fibers in this case), non homogeneous materials
are solved via \emph{``averaged equations''} so that the intricate
microstructure results in a continuum with averaged properties.  The
averaged equations and properties can be deduced in a rigorous way
through the theory of mixtures and averaging
operators~\citep{whitaker80}.

\subsection{Energy balance}

We will not enter in the details of all the derivations but only for
the averaged energy balance equation, which can be found in
Appendix~\ref{energy-bal}. The referred equation is 
\begin{equation} 
  \rho_s\Cp\dep Tt = \nabla \cdot ( k\nabla T- \rho_v\Vg (\Cpv T+\lambda+\Ql)) - 
      \dot m (\lambda+\Ql), \label{eq:balq} 
\end{equation}
where $T$ is temperature, $t$ time, $\Cp$ specific heat, $k$ thermal
conductivity, $\nabla$ the gradient operator, $\rho_s$ the density of
the dry board (solid phase), $\rho_v$ vapor density, $\Vg$ the volume averaged 
gas velocity, i.e. 
\begin{equation} \label{eq:vg} 
  \Vg=\epsilon\vg
\end{equation}
where $\vg$ is the velocity averaged on the phase, see
Appendix~\ref{energy-bal}), $\Cpv$ specific heat of vapor, $\dot m$
evaporation rate, $\lambda$ latent heat of vaporization of free water,
and $\Ql$ adsorption heat.  The main difference between this equation
and that presented by both \citet{carvalho98} and \citet{humphrey82}
is in the addition of the water evaporation heat term in the
convection term instead of considering the phase change effect only on
the temporal term. This term should be included because both phases,
the solid material and the vapor are in relative motion and we think
that its influence is not negligible in a high temperature process.

\subsection{Steam mass balance}

\citet{carvalho98} proposed the following steam
mass conservation equation
\begin{equation} 
    \dot m = \frac{\MMw}R \, \epsilon \nabla\cdot\LL[ - \Dv \nabla(\frac \Pv T) +
               \frac1\epsilon  \Vg \frac \Pv T \RR],
                  \label{eq:dotm0} 
\end{equation}
where $\MMw$ is the molecular weight of water, $R$ the gas constant,
$\epsilon$ board porosity, $\Dv$ diffusivity of water vapor in the
air/vapor mixture and $P_v$ vapor partial pressure. Considering
that the steam is treated as an ideal gas, then
\begin{equation} 
   \frac \Pv T = \MMw\muno R \rho_v,
\end{equation}
so it may be written, assuming $\epsilon$ constant, as
\begin{equation} 
   \dot m = \nabla\cdot\LL[ -\epsilon \Dv \nabla\rho_v + \Vg \rho_v \RR].
                  \label{eq:dotm} 
\end{equation}
This expression is preferable to  (\ref{eq:dotm0})
because it is written in a conservative form that is more agreeable
for a numerical treatment.  The left hand side term represents the
mass interfacial transport and those in the right hand side take into
account the mass diffusion and the mass convection. However, it should
be noted that this last expression does not have a temporal term 
as every consistent balance equation does.  For example, if evaporation is
not considered, then (\ref{eq:dotm}) is valid only for a steady
situation, which is not in general the case.  Then, we rewrite the
steam mass balance as
\begin{equation} 
   \epsilon \dep {\rho_v}t =  \nabla\LL[ \epsilon \Dv \nabla\rho_v 
             - \Vg \rho_v \RR] + \dot m.
                  \label{eq:dotm2} 
\end{equation}
This is another difference between our model and that proposed by
\citet{carvalho98}. 

\subsection{Gas mixture mass balance}

Finally, because the gas phase is composed of two main constituents,
steam and air, we may use an additional equation for the mass
transport of the whole gas phase. \citet{carvalho98} considered 
\begin{equation} 
   \dep Pt = -\frac 1\epsilon \nabla\cdot\LL( -\frac \KKg\mu \frac PT
                \nabla P \RR) \, T + \frac{\dot m}{\epsilon\MMa} TR +
                \frac PT \dep Tt, 
\end{equation}
where $P$ is the pressure of the gas phase and $\Kg$ the board permeability tensor. 
Again, assuming ideal gas law as the state equation for this phase,
\begin{equation} 
   \dep {}t \LL(\frac PT\RR) = -\frac{R}{\epsilon\MMa} \nabla\cdot\LL( -\frac
                \Kg\mu \rho_g \nabla P \RR)  + \frac{\dot m R}{\epsilon\MMa},
\end{equation}
we arrive to
\begin{equation} 
   \epsilon \dep {\rho_g} t = -\nabla\cdot (\rho_g \Vg) + \dot m.
           \label{eq:balg} 
\end{equation}

In order to close the system of equations we need to introduce a
relationship between $\dot m$ and $\deps \Pv t$. 
Consider the steam mass balance  (\ref{eq:dotm2})
and the relation
\begin{equation} 
  \rho_s \HH = \epsilon \rho_v + \rho_L,
\end{equation}
that represents the fact that the board moisture content $\HH$ is
composed of vapor and bound water $\rho_L$. If we assume that no
liquid phase is considered, then bound water may be transferred to the
gas phase only (solid to steam). So
\begin{equation} 
   \dot m = - \dep {\rho_L} t,
\end{equation}
and then
\begin{equation} 
\begin{split}
  \rho_s \dep {\HH}t &= \epsilon\dep{\rho_v} t + \dep{\rho_L} t \\
        &= \nabla\cdot\LL[ \epsilon \Dv \nabla\rho_v 
                - \Vg \rho_v \RR]. \label{eq:balh} 
\end{split}
\end{equation}
The air mass balance equation can be obtained by 
substracting (\ref{eq:dotm2}) from (\ref{eq:balg}) 
\begin{equation} 
   \epsilon \dep {\rho_a}t =  \nabla\cdot\LL[ -\epsilon \Dv \nabla\rho_v 
             - \Vg \rho_a \RR]. 
\end{equation}
Due to the fact that the mean macroscopic diffusive fluxes should be null 
\begin{equation} 
  \Dv \nabla\rho_v + \Da \nabla\rho_a = 0, 
\end{equation}
the air mass balance equation is transformed in the following expression
\begin{equation} 
   \epsilon \dep {\rho_a}t =  \nabla\cdot\LL[ \epsilon \Da \nabla\rho_a
             - \Vg \rho_a \RR], \label{eq:balai} 
\end{equation}
which is very similar to (\ref{eq:dotm2}) but here valid for the air.
Obviously, the  air transport equation has no evaporation term.

\section{Summary of equations and boundary conditions}\label{sec:summ}  

In order to clarify the mathematical model that is finally used for
the simulation of hot-pressing process we present the following brief
summary of partial differential equations.

\begin{itemize}
\item Energy balance equation 
\begin{equation} 
\rho_s\Cp\dep Tt = \nabla \cdot ( k\nabla T- \rho_v\Vg (\Cpv T+\lambda+\Ql)) - 
         \dot m (\lambda+\Ql).
\end{equation}

\item Water content balance equation
\begin{equation} 
  \rho_s \dep {\HH}t = \nabla\cdot\LL[ \epsilon \Dv \nabla\rho_v 
                - \Vg \rho_v \RR].
\end{equation}

\item Air mass balance equation
\begin{equation} 
   \epsilon \dep {\rho_a}t =  \nabla\cdot\LL[ \epsilon \Da \nabla\rho_a
             - \Vg \rho_a \RR]. 
\end{equation}
\end{itemize}

The boundary conditions are the following:
\begin{itemize}
\item At the press platen: 
\begin{equation} 
\begin{split}
   &T=\Tplaten(t),\ \ \textrm{air/water mixture in equilibrium with platen temperature}\\
   &\Vg\cdot\nor=0,\ \ \textrm{no mass flow across the platen}\\
   &\dep{\rho_a}z=0,\ \ \textrm{no air diffusion across the platen}\\
   &\dep{\rho_v}z=0,\ \ \textrm{no vapor diffusion across the platen}.
\end{split}
\end{equation}
\item At the center line ($r=0$), axial symmetry for all variables
\begin{equation} 
   \dep Tr = 0,\ \ \dep Hr = 0,\ \ \dep {\rho_a}r = 0.
\end{equation}
\item At the mid plane ($z=0$), symmetry for all variables 
\begin{equation} 
   \dep Tz = 0,\ \ \dep Hz = 0,\ \ \dep {\rho_a}z = 0.
\end{equation}
\item At the exit boundary ($r=R_{\textrm{ext}}$),
\begin{equation} 
\begin{split}
   \dep Tr &= 0,\ \ \textrm{null diffusive heat flux}\\
   P_v &= P_{v,\textrm{atm}},\ \ \textrm{equil. with
                   external air/water mixture,}\\
   P_a &= P_{a,\textrm{atm}},\ \ \textrm{equil. with
                   external air/water mixture.}
\end{split}
\end{equation}
\end{itemize}

\section{Numerical method}

The above system of equations contains three main unknowns, the
temperature, the moisture content and the air density representing the
dependent variables of the problem also called the state variable. In
this work we have used as independent variables the time and two
spatial coordinates (3D problems may be computed much in the same
way).  Due to the physical and geometrical inherent
complexity of this problem this may be computed only by numerical
methods.  For the spatial discretization we have employed finite
elements with multilinear elements for all the unknowns. Due to the
high convective effects the numerical scheme was stabilized with the
SUPG (for \emph{``Streamline Upwind - Petrov Galerkin''}) formulation
(see~\citet{BH_82}), otherwise spurious oscillations arise.  Once the
spatial discretization is performed the partial differential system of
equations is transformed into an ordinary differential system of
equations like
\begin{equation} 
   \dot U = R(U),
\label{eq:odesystem} 
\end{equation}
where $U$ is the state vector containing the three unknowns in each
node of the whole mesh. So, the system dimension is $3N$ where $N$ is
the number of nodes in the mesh.  The numerical procedure is as
follow: Knowing the state vector at the current time $(t^n)$,
i.e. $U_j(t^n) = [T_j,\HH_j,\rho_{a,j}](t^n)$ where $j$ represent an
specific node in the mesh. To get the residual, right hand side of
(\ref{eq:odesystem}) the following steps should be done:

{\raggedright
\begin{itemize}
\item Obtain air pressure from the gas state equation ($(T,\rho_a)\to P_a$).
\item Obtain relative humidity from sorption isotherms ($(T,\HH)\to\HR$).
\item Compute saturated vapor pressure from Kirchoff expression
                (\ref{eq:psat}) ($T->\Psat$).
\item Compute vapor pressure $P_v = \HR\Psat$. 
\item Obtain vapor density in air from vapor state equation $(T,\Pv)\to \rho_v$.
\item Compute coefficients from additional constitutive laws
        $(P,\HH,T)$ $\to$ $(D, \Kg, k_{x,y}, \Cp$).
\item Compute gradients of $T$, $P$ from nodal values at the 
        Gauss points using the  finite element interpolation.
\item Assemble the  element residual contributions in a global vector.
\end{itemize}
}

Once the whole residual vector at $t=t^n$ is computed the unknowns variables 
at the next time step is updated with
\begin{equation} 
   U^{n+1} = U^n + \Dt R(U^n).
\end{equation}
This kind of scheme, called {\sl explicit integration in time} is very
simple to be implemented but it has two major drawbacks, one is the
limitation of the time step to ensure numerical stability and the
other is the bad convergence rate for ill-conditioned system of
equations.  In this application the last disadvantage is very
restrictive because the characteristic times of each equation are very
different.  To circumvent this drawback we have implemented an
implicit numerical scheme
\begin{equation} 
   U^{n+1}- \Dt R(U^{n+1}) = U^n,
\end{equation}
where the residue is computed at the new state variable $U(t^{n+1})$
instead of using the current value $U(t^n)$. The non-linearities and
the time dependency of the state vector makes this implementation more
difficult and time consuming but more stable. This non-linear equation 
in $U^{n+1}$ is solved by the Newton method, which requires the
computation of the Jacobian
%
\begin{equation} 
    J = \dep RU.
\end{equation}
In order to avoid the explicit computation of the Jacobian, we compute
an approximate one by finite differences. 
\begin{equation} 
    J \approx J_{num} = \frac{R(U + \delta U) - R(U)}{\delta U}.
\end{equation}
This can be done element-wise, so the cost involved is proportional to
the number of elements in the mesh. In our application and despite of
the ill-condition of the problem we have found a good convergence at
each time step, in average 4 iterations per time step to reduce 10
orders of magnitude in the global residue.

\section{Physical and transport properties}

\subsection{Thermal conductivity}

Following \cite{humphrey82} the influence of board
density, moisture content and temperature on the thermal conductivity
is considered as an independent correction factor obtained experimentally
\begin{equation}
\begin{split}
     \kappa_z &= F_{\kappa,H}\, F_{\kappa,T} 
               \, (\nexp{1.172}{-2} + \nexp{1.319}{-4} \, \rho_s), \\
        F_{\kappa,H} &= 1+ \nexp{9.77}{-3} \, (H-12),\\
        F_{\kappa,T} &= (T-20)\times 1.077 \times 10^{3} + 1,
\end{split}
\end{equation}
where $\kappa_z$ is the thermal conductivity in the pressing direction in 
$\rm W/m \gK$ and $\rho_s$ is the oven dry density of the material in $\rm Kg/m^3$.
The moisture correction factor $F_{\kappa,H}$ (\citet{kollmann56} and
\citet{humphrey82}) assumes $H$, the moisture content of the board
material, in $\%$, and the temperature correction factor
$F_{\kappa,T}$ (\citet{kuhlmann62} and \citet{humphrey82}) assumes $T$
in $\gC$. 

\subsubsection{Heat flux direction correction}
By far the greater part of conductive heat translation takes place in
the vertical plane. However the energy lost from the mattress is
largely the result of radial vapor migration from the center toward
the atmosphere. The associated horizontal relative humidity gradient
lead to a horizontal temperature gradient. Even though this gradient
is always lower than the vertical one its influence should be taken
into account if multidimensional analysis is required. 
\citet{ward63} made experimental measurements and they observed that
at a first glance a factor of approximately 1.5 may be a good initial
guess before doing some extra measurements. Then
\begin{equation} 
        \kappa_{xy} = 1.5 \kappa_z .
\end{equation}

\subsection{Permeability}

The evaporation and condensation of water changes the vapor density
and consequently its partial pressure in the voids within the
composite. A vertical pressure gradient leads to the flow of water
vapor from the press platens toward the central plane of the
board. At the same time an horizontal vapor flow is set up in
response to the pressure gradient established in the same
direction. The relation between the pressure gradient and the flow
features may be assigned to the material permeability.  Permeability
is a measure of the ease with which a fluid may flow through a porous
medium under the influence of a given pressure gradient. Different
mechanisms may be involved in this flow, a viscous laminar flow, a
turbulent flow and a slip or Knudsen flow. In this study only the
first type is included with the assumption that Darcy's law is
obeyed. This may be written as
\begin{equation} 
        \Vg = -\frac{!K_g}{\mu_g} \nabla p
\end{equation}
where $\nabla p$ is the driving force and $\mathbf{v}_{g}$ is the flow
variable. ${!K_{g}/\mu_{g}} $ is called the superficial
permeability with $\mu_g$ the dynamic viscosity of the gas phase and
$!K_{g}$ the specific permeability.  

The kinetic theory of gases suggests that at
normal pressure the viscosity is independent of pressure and it varies
as the square root of the absolute temperature,
\begin{equation} 
        \mu \propto \sqrt{T} \ \ \ \ \   10^3 < p < 10^6.
\end{equation}

Corrections with the absolute temperature are often considered by Sutherland law
\begin{equation} 
        \mu = \frac{B \,\, T^{3/2}}{T+C},
\end{equation}
with $B$ and $C$ characteristic constants of the gas or vapor with
$\mu$ in $Kg \, m^{-1} \, s^{-1}$. Values for $B$ and $C$ are
available from \citet{keenan66}. For this application the
following expression was adopted,
\begin{equation}
    \mu = 1.112 \times 10^{-5} \times \frac{(T+273.15)^{1.5}}{(T+3211.0)},
\end{equation}
with $T$ in $\gC$.

\subsubsection{Variation of vertical permeability with board material density}

Even though we consider that the press is closed and the density
profile is set up and fixed it is included here some conclusions from
\cite{humphrey82} with results obtained by \citet{denisov75} for 19 mm
boards and others from \citet{sokunbi78}. The data to be fitted are the
following

\begin{table}[htb]
\centering
\begin{tabular}{|c|c|}
\hline
Mean density     & Mean vertical \\
 $[\rm Kg/m^3]$  &  permeability\\ 
                 &  $\rm [m^2 \times 10^{15}] $\\
\hline
\hline
425   &  64\\
\hline
475   &  40\\
\hline
525   &  24\\
\hline
575   &  16\\
\hline
625   &  11\\
\hline
675   &  7 \\
\hline
725   &  5 \\
\hline
775   &  3 \\
\hline
825   &  2 \\
\hline
875   &  2 \\
\hline
\end{tabular}
\caption{Table I: Vertical permeability density correction data}
\end{table}

\subsubsection{Horizontal permeability }

Sokunbi measures included in figure 2.7 of \cite{humphrey82} shows the
relation between the board thickness in mm with horizontal
permeability.  For approximately $15 \rm mm$ board thickness and
$\rho_s = 586 \rm Kg/m^3$ the horizontal permeability is 59 times the
vertical value, in agreement with the values assumed by
\citet{carvalho98}.

\subsection{Steam in air diffusivity}

The interdiffusion coefficient of steam in air can be calculated from
the following semi-empirical equation \citep{stanish86},
\begin{equation} 
  D_a = \nexp{2.20}{-5} \, \LL(\frac{101325}P\RR) \, \LL(\frac
                    T{273.15}\RR)
\end{equation}
where the diffusivity is in $\rm m^2/sec$, pressure in $\rm N/m^2$ and $T$ in
$\gK$. 

\subsection{Vapor Density}

For the pressure range likely to occur during hot pressing (between
$10^3$ and $3\times 10^5\rm N/m^2$ ) a linear relationship between saturated
vapor pressure and vapor density may be assumed. Fitting
experimental data \citet{humphrey82} proposed the following expression
\begin{equation} \label{eq:H_2_18}
  \rho_v = P_{sat} \,\nexp{6.0}{-8} \, HR,
\end{equation} 
with $\rho_v$ in $\rm Kg/m^3$, $P_{sat}$ in $\rm N/m^2$ and the relative
humidity $HR$ in \%. This can be deduced from the relative humidity
definition 
\begin{equation} 
        HR = P_v / P_{sat}, 
\end{equation}
and applying ideal gas law for gaseous phase 
\begin{equation} 
        \frac{P_v}{\rho_v } = \bar R/\MMw T .
\end{equation}
Taking $ \bar R = 8314 \rm J/Kmol/K$ and $\MMw= 18 \rm Kg/Kmol$ with $T
\approx 360 \gK$ we obtain (\ref{eq:H_2_18}). 

\subsection{Saturated vapor pressure}

Following the Kirchoff expression with data presented by
\citet{keenan66} we include here the following equation,
\begin{equation} \label{eq:psat}  
        \log_{10} P_{sat} = 10.745 - (2141.0/(T+273.15)),
\end{equation}
with $P_{sat}$ in $\rm N/m^2$ and $T$ in $\gC$.

\subsection{Latent heat of evaporation and heat of wetting}

Using Clausius-Clapeyron equation in differential form and after some
simplification the latent heat of vaporization of free water may be
written as
\begin{equation} 
   \lambda = \nexp{2.511}6 - \nexp{2.48}3 \, T ,
\end{equation}
with $\lambda$ in $\rm J/Kg$ and $T$ in $\gC$.

For the differential heat of sorption we follow \citet{humphrey82} that used
(see \citet{bramhall79})
\begin{equation} 
  \Ql = \nexp{1.176}6 \, \expe{-0.15 H},
\end{equation}
with $Q_l$ in $\rm J/Kg$ and $H$ in \%.

\subsection{Specific heat of mattress material}

It is computed by adding the specific heat of dry wood and that of
water according to the material porosity and assuming full saturation.
The specific heat of dry mattress material is taken as $1357 \rm J/Kg/K$
and the specific heat of water has been taken to be $4190
\rm J/Kg/K$. From \citet{siau84} the expression for specific heat of moist %
wood is
\begin{equation} 
  C_p = 4180 \frac{ 0.268 + 0.0011 (T-273.15) + H}{1+H},
\end{equation}
$T$ in $\gK$ and $H$ in \%. 

\subsection{Porosity}

According to \citet{humphrey82} the volume of voids within the region may be
computed with
\begin{equation} 
        \epsilon = \frac{V_{\mathrm{voids}}}V  = (1 - \frac{\rho}{\rho_s} ),
\end{equation}
where $\epsilon$ is the porosity, $\rho$ is the density of the region
and $\rho_s$ is the dry density of the board material.

In \citet{carvalho98} they included the expression from \citet{suzuki89}
\begin{equation} 
  \epsilon = 1 - \rho_s \frac{1/\rho_f + y_r/\rho_r}{1+ y_r},
\end{equation}
where $\rho_r$ is the cure resin density, $\rho_f$ is the oven dry
fiber density and $y_r$ is the resin content (resin weight/board
weight).  In \citet{carvalho98}, the authors used $y_r= 8.5\%$ with $\rho_f=900
\rm Kg/m^3$ and $\rho_r=1100 \rm Kg/m^3$.

\section{Numerical results}

\begin{figure}
\centerline{\includegraphics{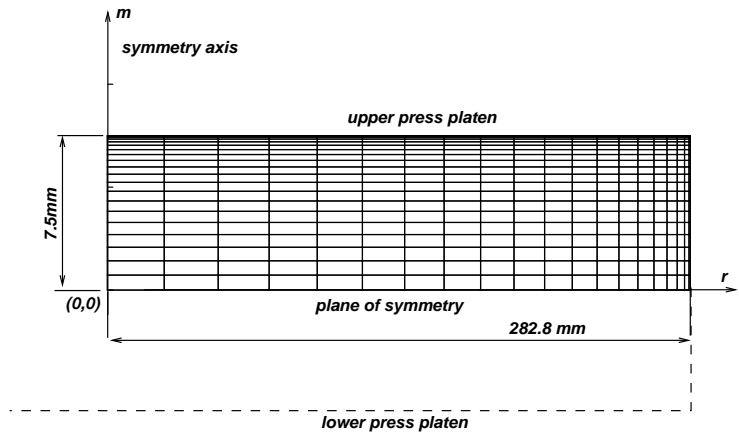}}
\caption{Finite element mesh}
\label{fg:mesh}
\end{figure}

In this section we present some results that can be compared with
those reported in \citet{humphrey82}.  This numerical experiment
allows the validation of the mathematical model and its numerical
implementation for future applications to hot-pressing process
simulation and control.  This experiment consists of a round
fiberboard of 15 mm of thickness and 0.2828 m of radius that according
to its axisymmetrical geometry needs as spatial coordinates only the
radius $r$ and the axial coordinate $z$ assuming no variation in the
circumferential coordinate. The axisymmetrical domain is discretized
in $20 \times 20$ elements in each direction with a grading toward the
press platen and the external radius as may be visualized in
figure~\ref{fg:mesh}. In order to follow the same assumptions as in
that work we fixed the air density to a very low value
$(\rho_a \approx \nexpE{-6})$, as if the press was close with no air
inside the fiberboard. The boundary conditions are as in
section~\ref{sec:summ}. 

For the press platen temperature we have applied a ramp from $30
\gC$ at $t=0$ to $160 \gC$ at $t=72 $ seconds with a least
square fitting from data in \citet{humphrey82}. 
As initial conditions we have assumed a uniform temperature of
$T(t=0)=30 \gC$ in the whole solid material with a uniform moisture
content of $H=11 \%$. The external atmosphere was considered to be at
$T_\atm=30\gC$, $\HR_\atm=65\%$, in such a way that the internal
moisture content is, at the initial state, in equilibrium with the
external atmosphere.

In the next section we include the results obtained by using the model
above cited to the original problem presented in \citet{humphrey82}. Next we
present some other experiments showing some phenomena that deserve
more attention for futue studies.

\subsection{Original numerical experiment}

\begin{figure}
\centerline{\includegraphics{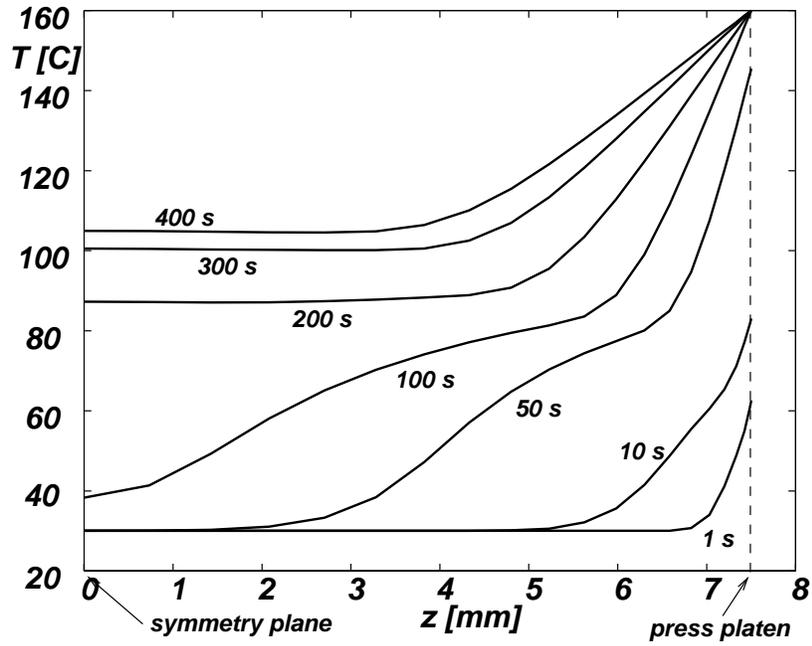}}
\caption{Axial temperature profile at $r=0$}
\label{fg:T_zt_r0}
\end{figure}

\begin{figure}
\centerline{\includegraphics{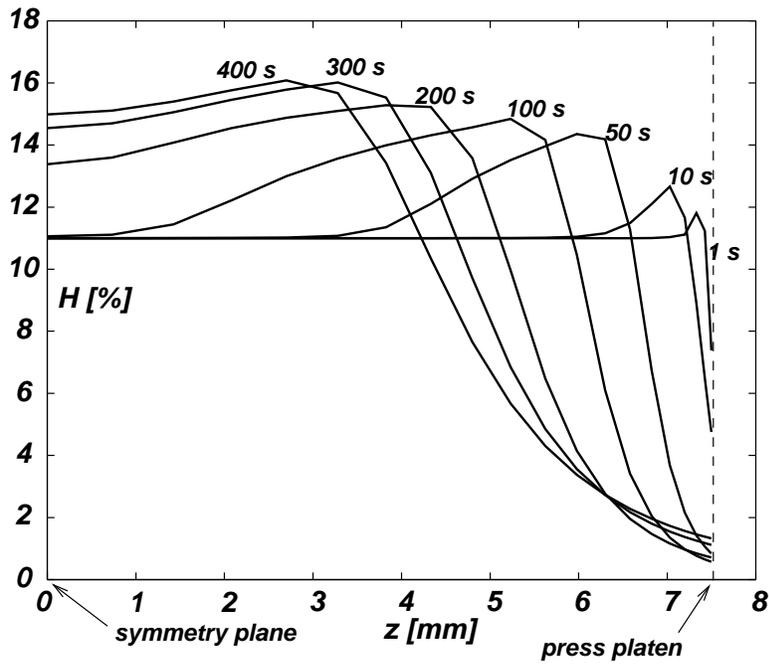}}
\caption{Axial moisture content profile at $r=0$}
\label{fg:H_zt_r0}
\end{figure}

\begin{figure}
\centerline{\includegraphics{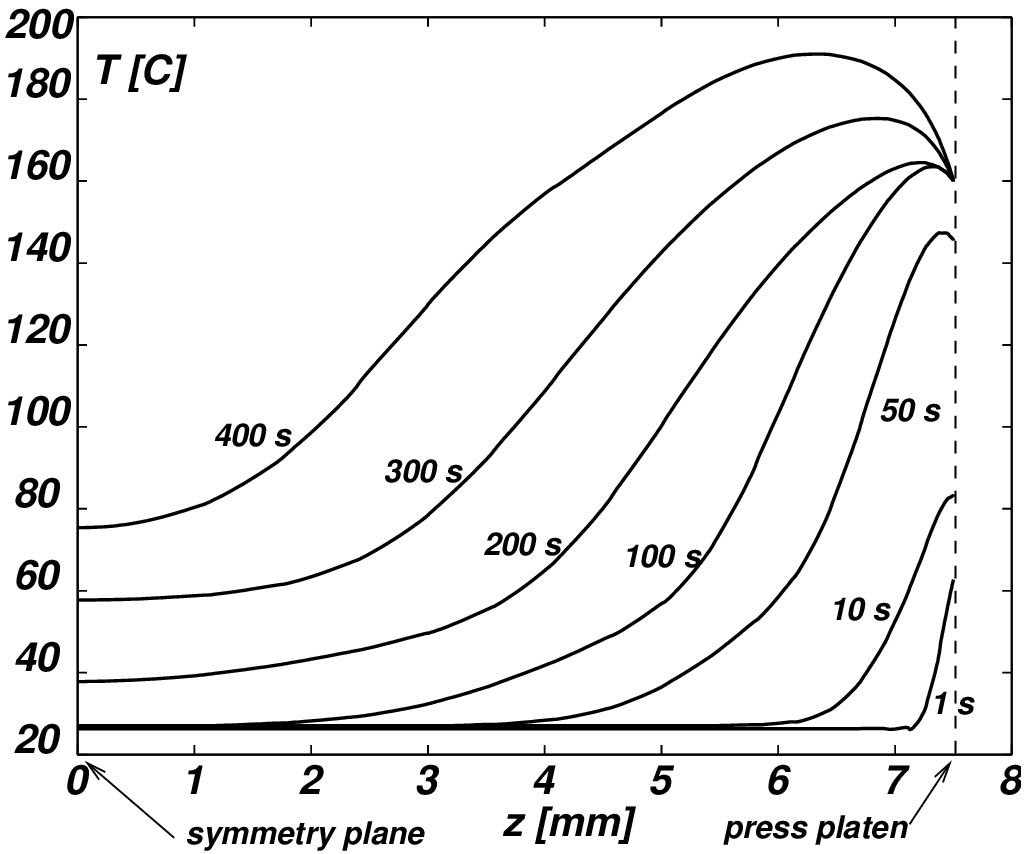}}
\caption{Axial temperature profile at $r=R=0.2828 \rm m$}
\label{fg:T_zt_rR}
\end{figure}

\begin{figure}
\centerline{\includegraphics{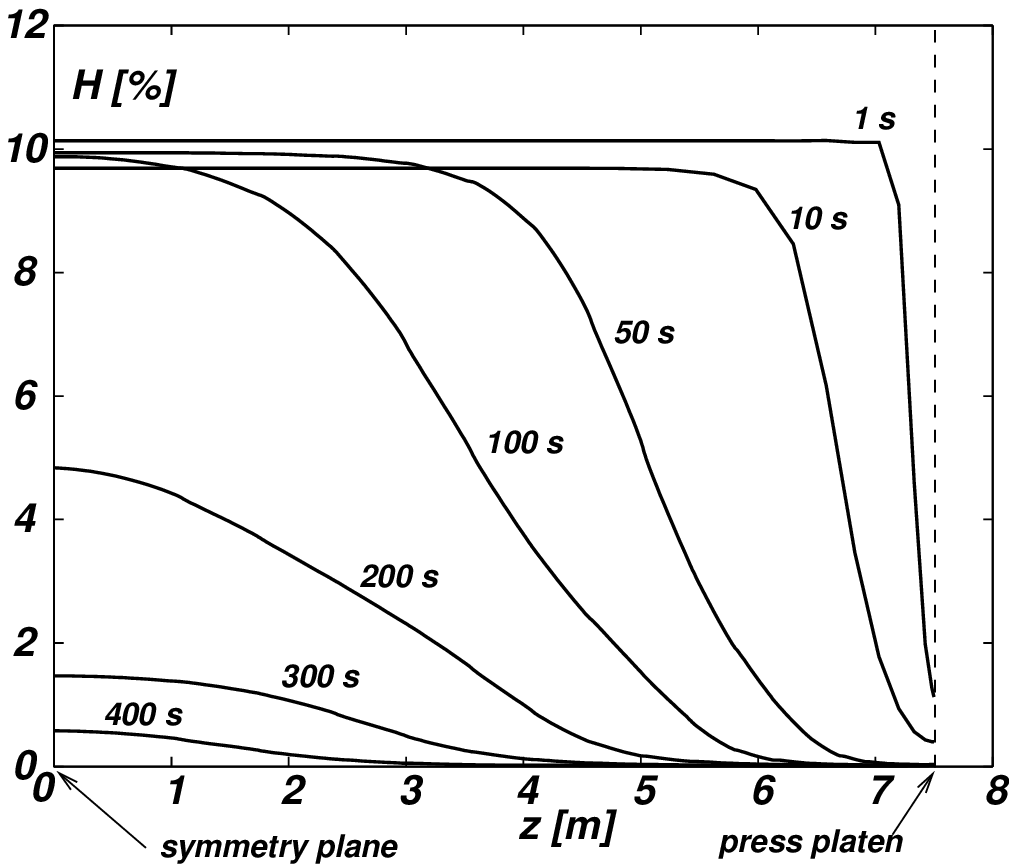}}
\caption{Axial moisture content profile at $r=R=0.2828 \rm m$}
\label{fg:H_zt_rR}
\end{figure}

\begin{figure}
\centerline{\includegraphics{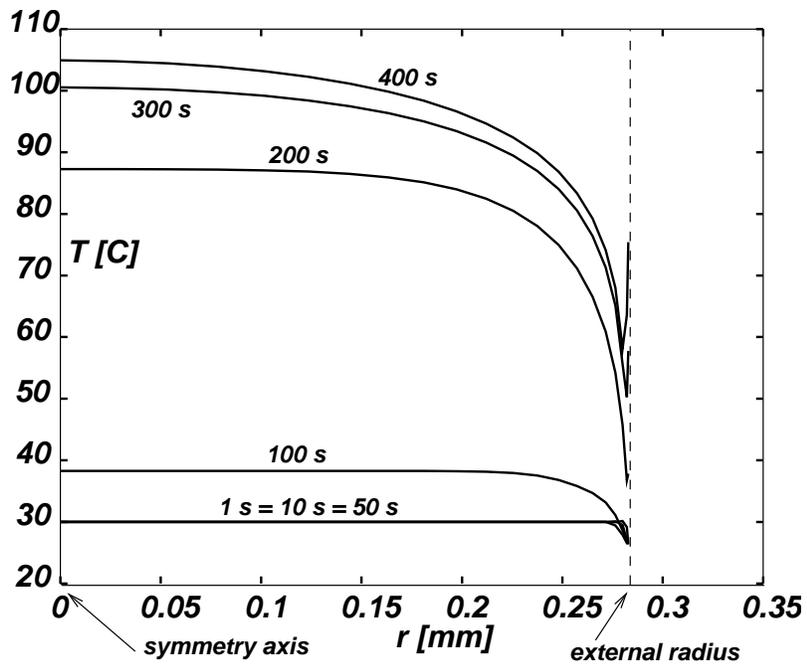}}
\caption{Radial temperature profile at $z=0$}
\label{fg:T_rt_z0}
\end{figure}

\begin{figure}
\centerline{\includegraphics{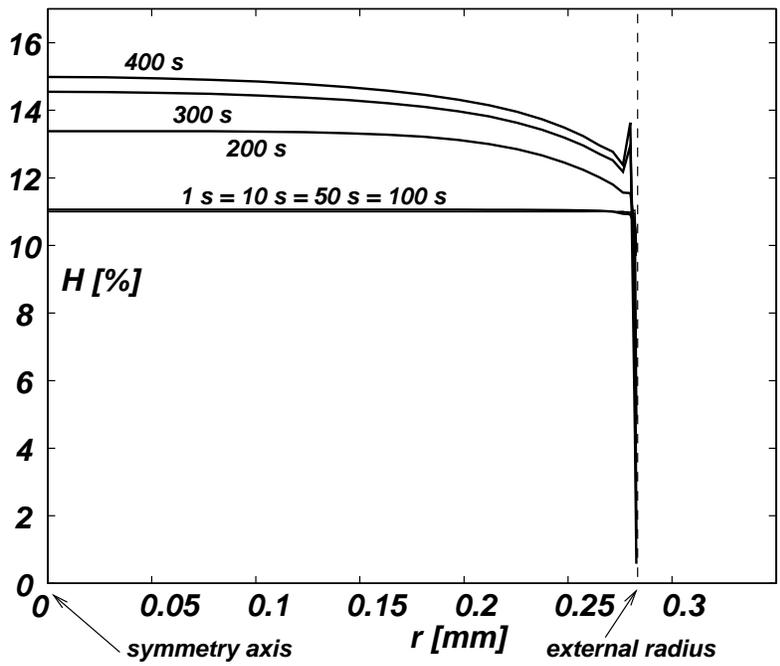}}
\caption{Radial moisture content profile at $z=0$}
\label{fg:H_rt_z0}
\end{figure}

\begin{figure}
\centerline{\includegraphics{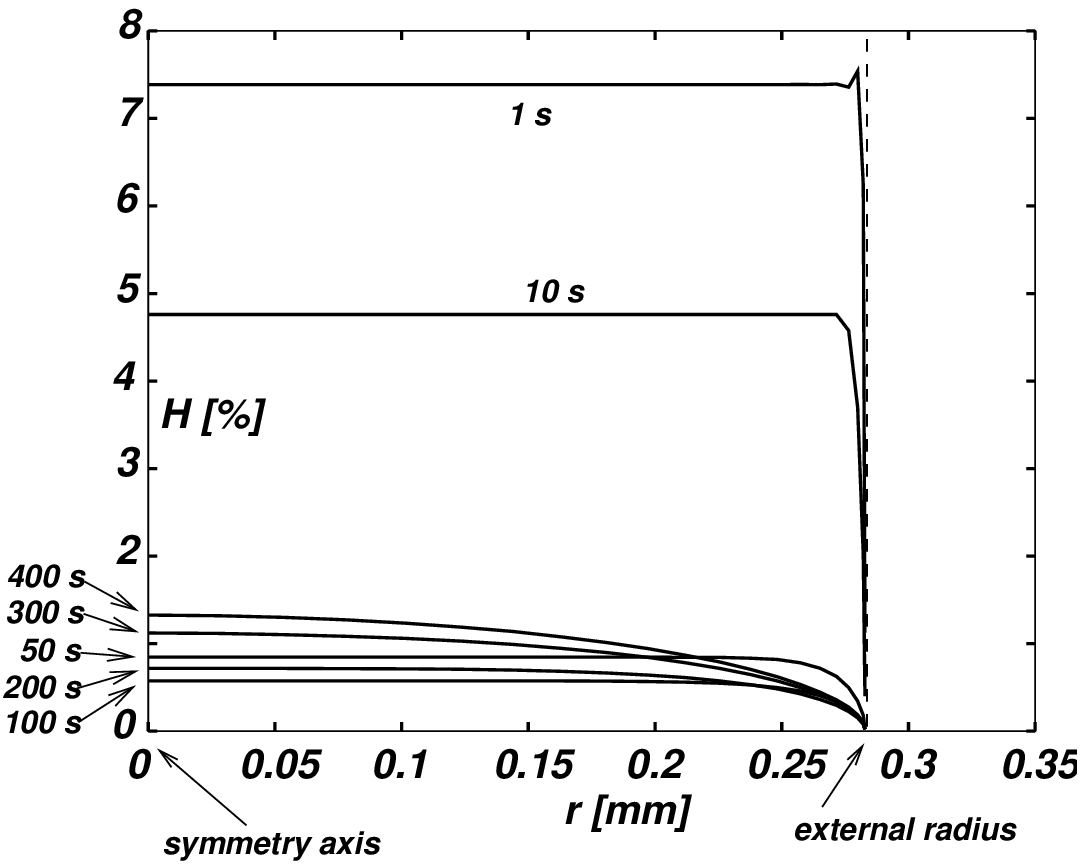}}
\caption{Radial moisture content profile at press platen}
\label{fg:H_rt_zP}
\end{figure}

Figure~\ref{fg:T_zt_r0} shows the temperature distribution in time and
with the axial coordinate at $r=0$ (centerline). We can note the
penetration in the axial direction of the temperature profile in time
for $t=1,10,50,100,200,300$ and $400$ seconds.
Figure~\ref{fg:H_zt_r0} shows the same kind of plot for moisture
content. For $r$ not to close to the external radius, the problem is
almost one dimensional in the $z$ direction. As the thermal front
penetrates into the board water evaporates. This vapor advances to
lower pressure regions near the symmetry plane and, as it encounters
lower temperatures, it condenses releasing heat. This process can be
clearly seen from the wave in moisture content
(figure~\ref{fg:H_zt_r0}) exceeding the initial water content of
11\% and results in an improvement in the heat transfer with respect
to the pure conduction case. 
Also this phenomenon is responsible of the change in curvature of thee
temperature curves, mainly at $t=10$, 50, and 100~sec 
(see figure~\ref{fg:T_zt_r0}). 
The total water content in the board at a given instant can be
found by integrating the bound water content and the water in vapor
phase. However, this last is negligible. We can see in
figure~\ref{fg:H_zt_r0} that the depression in water content near the
board (for instance at $t=400\rm sec$), is larger that the water
enrichment in near the center plane. This is due to water migration
from the center of the board to the external radius, where it flows to
the external atmosphere. 
The following figures show similar plots but at different locations,
\begin{itemize}
\item Figure \ref{fg:T_zt_rR}: temperature at external radius 
\item Figure \ref{fg:H_zt_rR}: moisture content at external radius 
\item Figure \ref{fg:T_rt_z0}: temperature at axial centerline $(z=0)$ 
\item Figure \ref{fg:H_rt_z0}: moisture content at axial centerline $(z=0)$
\item Figure \ref{fg:H_rt_zP}: moisture content at press platen
\end{itemize}

\begin{figure}
\centerline{\includegraphics{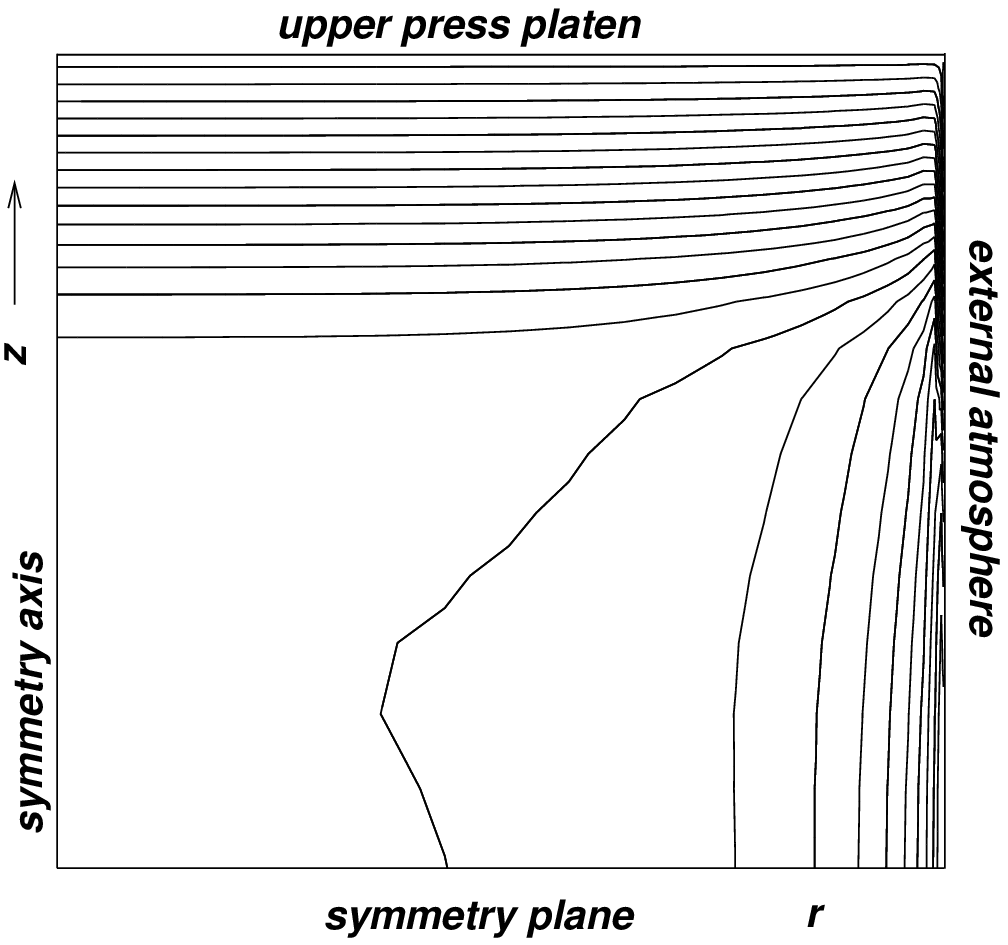}}
\caption{Isotherms at t=200 seconds ($z$ and $r$ axis not to scale.)}
\label{fg:T_iso_t200}
\end{figure}

Figure \ref{fg:T_iso_t200} shows several isotherms at $t=200$ seconds
distributed in the $r,z$ plane.

\begin{figure}
\centerline{\includegraphics{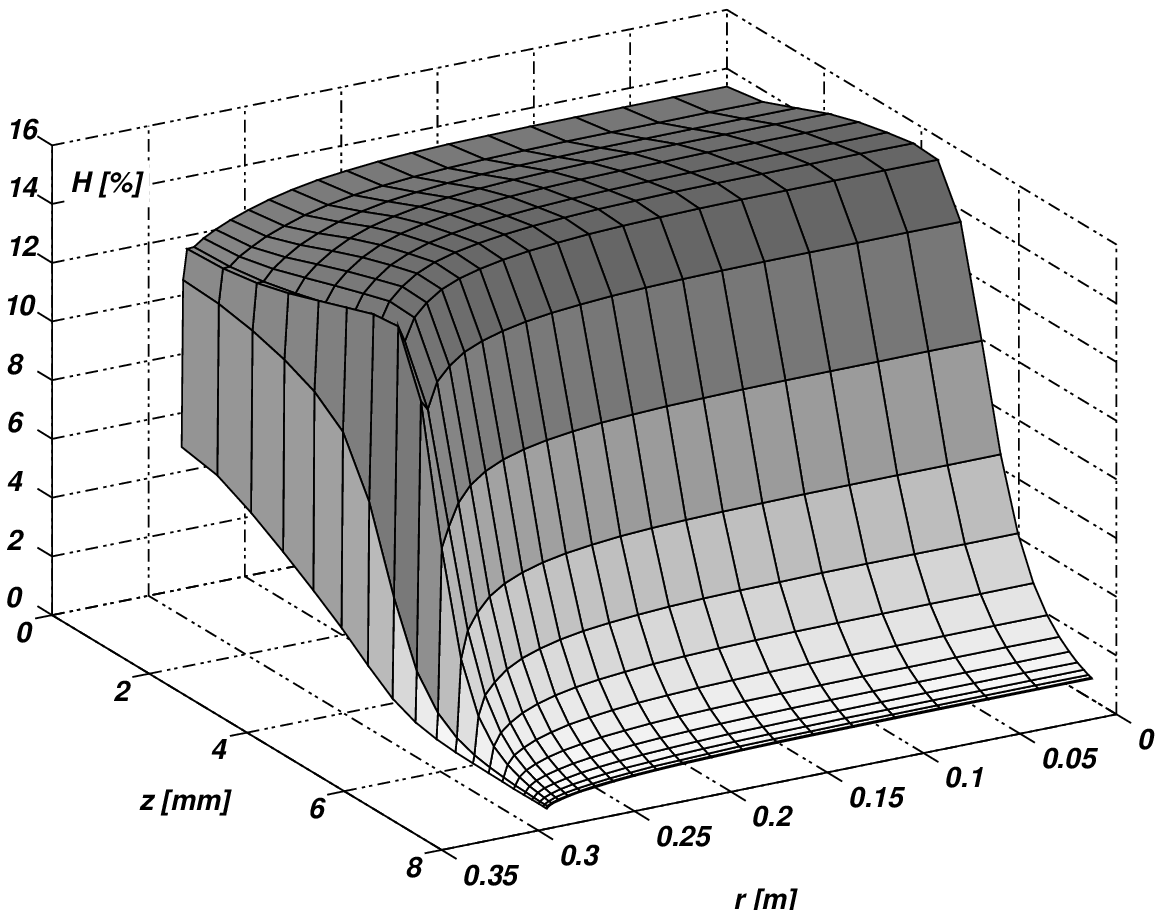}}
\caption{3D view for moisture content at t=200 seconds}
\label{fg:H_3D_t200}
\end{figure}

Figure \ref{fg:H_3D_t200} shows a three-dimensional view for moisture
content represented by the third coordinate axis at $t=200$ seconds as
a function of $r,z$.

These results are in good agreement with those presented by
\citet{humphrey82}.  However the cited author did not present his
results at some locations that in our opinion should be treated with
some care, for example at the external radius.

\subsection{Further numerical experiments}

In \citet{humphrey82}, results for moisture and temperature in the
vertical and radial directions at both central planes, $r=0$ and $z=0$
respectively are included. No mention about the vertical distribution at
$r=R=0.2828 \rm m$ or about the moisture content at press
platen. Moreover, he had used a uniform mesh of $10 \times 10$
elements without showing what happen at the last annuli of elements
corresponding to the external radius.

\begin{figure*}[htb]
\centerline{\includegraphics{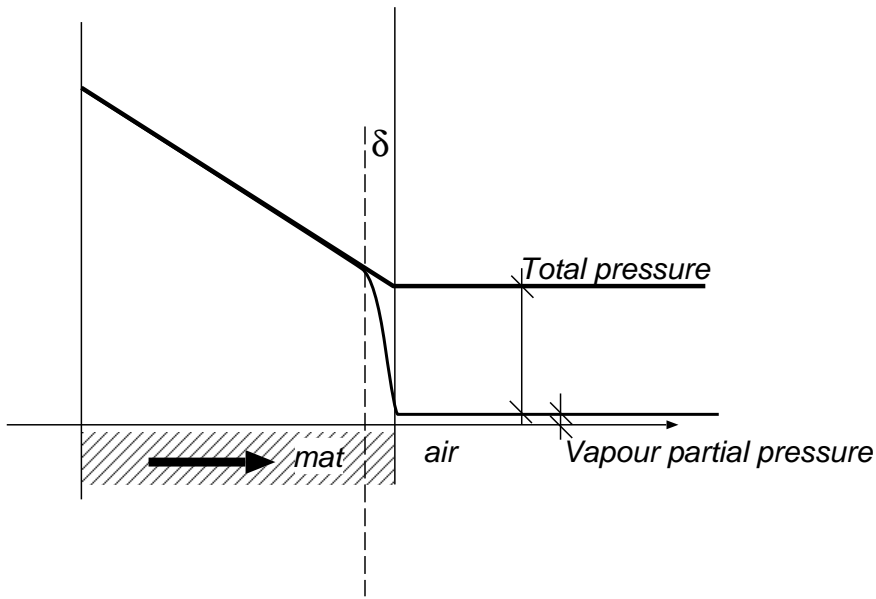}}
\caption{Boundary layer in vapor partial pressure at external radius}
\label{fg:blay}
\end{figure*}

Our results present some overshooting in the temperature profile very
close to the radial exit contour and at the first moment we though
about a spurious numerical problem, but it is due to large variations
of the magnitude of vapor pressure and density at the boundary.  In
typical runs, vapor pressure varies from near 2~atm in the center of
the board to 0.01~atm at the external radius. We think that this
problem will be fixed if we solve for the air density also, but then a
very fine grid will be necessary at the exit boundary, since large
variations of the vapor molar fraction are expected. (see
figure~\ref{fg:blay}). Molar fraction varies from nearly 1 at the
interior of the mat to a 2\% at the external atmosphere. This
variation is produced in a thin layer of width $\delta$ proportional
to the diffusivity of vapor in air which is very small.

\section{Conclusion} 

We presented a numerical model for the heat and mass transfer in the
hot-pressing model of a MDF fiberboard. The model includes convective
effects on the phase change term and also a conservative numerical
treatment of the resulting system of partial differential
equations. Convective effects are responsible of an increase in heat
transfer from the platen to the center of the board due to water
vapor evaporation and condensation. Two-dimensional simulations allow
to estimate border effects. 

\subsection*{ACKNOWLEDGEMENTS}

This work has received financial support from Consejo Nacional de
Investigaciones Cient\'\i ficas y T\'ecnicas (CONICET, Argentina,
grants PIP 0198/98, PIP 0266/98), Universidad Nacional del Litoral
(Argentina, grants CAI+D 96-004-024, CAI+D 2000/43) and ANPCyT
(Argentina, grants PICT 12-0051, PICT 12-06973, PID-99/74,
PID-99/76).
We made extensive use of freely distributed software as GNU/Linux~OS,
GNU/Octave and others.

\appendix
\section{List of symbols}

\subsection{Physical constants and quantities} 

\begin{itemize}
\item [$D$:] diffusivity  \([\rm m^2/s]\)
\item [$\Dt$:] time step  \([\rm s]\)
\item [$\epsilon$:] material porosity (dimensionless)
\item [$F_k$:] correction factor for thermal conductivity (dimensionless)
\item [$\mathrm{Fo}$:] Fourier number (dimensionless)
\item [$h$:] enthalpy
\item [$\HH$:] water content in \% weight 
\item [$\HR$:] relative humidity [\%]
\item [$J$:] Jacobian matrix
\item [$k$:] thermal conductivity \([\rm W/m\gK]\)
\item [$!K_g$:] specific permeability  \([\rm m^2]\)
\item [$\lambda$:] latent heat of vaporization of free water \([\rm J/Kg]\)
\item [$\dot m$:] evaporation rate \([\rm Kg/m^3s]\) 
\item [$\textrm{MM}$:] air molecular weight [$\rm Kg/Kmol$]
\item [$\mu$:] dynamic viscosity  \([\rm Kg/m s]\)
\item [$N$:] number of nodes in the finite element mesh
\item [$P$:] pressure \([\rm N/m^2]\)
\item [$\Ql$:] adsorption heat \([\rm J/Kg]\)
\item [$R$:] gas constant for air (\(8314 \rm \rm KJ/Kmol\gK\)
\item [$R$:] vector of residuals for the discrete model 
\item [$r$:] radial coordinate [mm]
\item [$\rho$:] density \([\rm Kg/m^3]\)
\item [$T$:] temperature \([\rm \gK]\)
\item [$U$:] state vector for the discrete model
\item [$V$:] volume \([\rm m^3]\)
\item [$!v$:] gas mixture velocity  \([\rm m/s]\)
\item [$V_{\mathrm{voids}}$:] volume of voids \([\rm m^3]\)
\item [$y_r$:] resin content [weight \%]
\item [$z$:] coordinate normal to the plate [mm]

\end{itemize}

\subsection{Indices} 

\begin{itemize}
\item [eff:] effective quantities (averaged for the gas/solid mixture)
\item [$f$:]  fiber
\item [$g$:]  gas phase, (air/water mixture)
\item [$j$:]  nodal index
\item [$L,l$:] bound water
\item [platen:] quantity evaluated at the press platen
\item [$r$:]  resin
\item [ref:] reference state
\item [$s$:] solid phase
\item [sat:] saturated atmosphere 
\item [$v$:] water vapor
\item [w:] water

\end{itemize}

\subsection{Mathematical symbols} 

\begin{itemize}
\item [$\nabla$:]  gradient (nabla) operator
\item [$\dot{(\ )}$:] temporal derivative
\end{itemize}

\section{Derivation of the averaged energy balance equation}
\label{energy-bal}

The microscopic energy balance equation in the gas phase is 
\begin{equation} \label{eq:enermicro} 
  \dep{}t (\rho h) + \nabla\cdot(\rho\B v h) 
        = -\nabla\cdot (k\nabla  T).   
\end{equation}
where $h$ is enthalpy. 
For the other phases (solid and bound water) a similar expression
holds, but neglecting the advective term. Applying the volume average
operator~\citep{whitaker80}, we arrive to the following equation
averaged on the gas phase
\begin{equation} 
  \dep{}t (\eg \av{\rho h}{g})+\nabla\cdot (\eg \av{\rho h \B v}g) =
         \nabla\cdot(\eg \av{k \nabla T}g) + Q_g,
\end{equation}
$\eg$ is the volumetric fraction of phase $g$ (i.e. gas) and $\av Xg$
is the average of quantity $X$ on the volume occupied by phase $g$
\begin{equation} 
  \av Xg = \frac1{\Omega_g} \int_{\Omega_g} X \di\Omega.
\end{equation}
The term $Q_g$ is the total enthalpy flux through the solid-gas
interface $\Gamma$
\begin{equation} 
   Q_g = \int_{\Gamma} (\rho h)_g (\B v-\B w)\cdot\nor \di\Gamma,
\end{equation}
where $(X)_g$ is the value of property $X$ on the $g$ side of the
interface and $\B w$ is the velocity of the interface. Assuming that
$h_g$ is constant on all $\Omega_g$ (for a certain volume control) then 
\begin{equation} 
\begin{split}
   Q_g &= \av{h}g \int_{\Gamma} (\rho )_g (\B v-\B w)\cdot\nor
                             \di\Gamma,\\
       &= \av{h}g \dot m,
\end{split}
\end{equation}
where $\dot m$ is the rate of mass of water being evaporated. A common
drawback of averaged equations is that, when products of variables
like $\rho h$ appear in the microscopic equation, the average of the
product $\av{\rho h}g$ is obtained in the averaged equation. Now, it
is not true that 
\begin{equation} \label{eq:correl} 
  \av{\rho h}g =   \av{\rho}g \, \av hg, 
\end{equation}
so that the averaged equation contains more unknowns than the
original equation. A common assumption is that no
correlation exists between variables and so (\ref{eq:correl}) is
approximmately valid. This can be justified, for instance, if the
variations of each quantity around the mean is small. 

Then, applying the volume average operator over the gas, solid, and
bound water phases and assuming no correlations between variables we
obtain the following averaged equations for the three phases
\begin{equation} 
\begin{split}
  \dep{}t (\eg \rhog h\G)+ \nabla\cdot(\eg\rhog\vg h\G) &=
         \nabla\cdot(\eg\av{k\nabla T}g) - \dot m h\G  \ \ \textrm{gas phase},\\
  \dep{}t (\es \rhos h\So)&= 
         \nabla\cdot(\es\av{k\nabla T}s)  \ \ \textrm{solid phase and
         },\\
  \dep{}t (\el \rhol h\Li)&= \dot m \, h\Li  \ \ \textrm{liquid phase.}
\end{split}
\end{equation}
Also, the average operator is dropped from here on, and a subindex $g$
or $s$ implies averaging on that phase. Also, $\Vg$ the volume averaged 
gas velocity is 
\begin{equation} 
  \Vg = \frac1{\Omega} \int_{\Omega_g} !v \di\Omega = \eg\vg. 
\end{equation}
%
Note that in the body of the text $\Vg$ is used instead of $vg$. 
Now, summation of these three equations gives
\begin{multline}
  \dep{}t (\eg\rhog h\G+\es\rhos h\So+\el\rhol h\Li) +
       \nabla\cdot(\eg\vg\rhog h\G) = \\
           \nabla(\eg\av{k\nabla T}g+\es\av{k\nabla T}s) 
             + \dot m \, (\hl-\hg).
\end{multline}
Now, 
\begin{equation} 
  \eg\av{k\nabla T}g+\es\av{k\nabla T}s = \keff \nabla \ave{T},
\end{equation}
where $\keff$ is the average conductivity of the solid+water+gas
mixture.  The gas is assumed as an ideal mixture, so that the enthalpy
is the sum of the enthalpy of its constituents, and neglect the
contribution of the air constituent so that
\begin{equation} 
  \rho_g \hg = \rho_a h_a + \rho_v h_v \approx \rho_v h_v.
\end{equation}
Taking a reference state for the entalphy at a point on the adsorbed
state
\begin{equation} 
   h_v  = \Cpv (T-\Tref) + \lambda+Q_l.
\end{equation}
We also neglect the entalphy of the gas phase with respect to the
solid+water phases and put
\begin{equation} 
   \es\rhos h\So + \el\rhol h\Li = 
                  (\rho\Cp)_\effe \, (T-\Tref) ,
\end{equation}
where $\rho_\effe$, and $\Cp_\effe$ are averaged properties for the
moist board, as a function of temperature and moisture content. 
Finally, the averaged equation is
\begin{multline}
   \dep{}t [(\rho\Cp)_\effe \, (T-\Tref)] 
         + \nabla\cdot[\eg\rho_v\vg \,(\Cpv(T-\Tref)+\lambda+Q_l)] = \\
              \dive(\keff T) - \dot m (\lambda+Q_l)
\end{multline}
This is equivalent to (\ref{eq:balq}) through relation (\ref{eq:vg})
and assuming $\Tref=0\gC$. 

\bibliographystyle{plainnat}
\bibliography{hotpress}
\end{document}